\documentclass[final,leqno]{siamltex}

\usepackage{ifthen}

\usepackage{amsmath,amsfonts,amssymb,algorithm,algorithmic,psfrag}
\usepackage{array}
\usepackage{psfrag}
\usepackage{subcaption}
\usepackage{bm}
\usepackage{color}
\usepackage{enumitem}
\usepackage{graphicx}

\usepackage{xypic}

\newcommand{\A} {{\mathbb A}}
\newcommand{\B} {{\mathbb B}}

\newcommand{\R} {\mathbb R}
\newcommand{\M} {\mathbb M}
\newcommand{\mesh} {\mathcal M}
\newcommand{\elt}{\mathfrak m}
\newcommand{\face}{\mathfrak f}

\newcommand{\Div}{\textnormal{div\,}}
\newcommand{\Curl}{\textnormal{curl\,}}

\newcommand{\p}{\partial}

\newcommand{\cancel}[1]{}
\newcommand{\raw}{\rightarrow}

\newcommand{\veps}{\varepsilon}

\newcommand{\hdiv}{H(\textnormal{div})}

\newcommand{\hcurl}{H(\textnormal{curl})}

\newcommand{\sign}{\text{sign }}

\newcommand{\ds}{\displaystyle}

\newcommand{\whitC}[3]{\lambda_{#1}\nabla \lambda_{#2}\times\nabla \lambda_{#3}}

\newcommand{\whit}{\mathcal{W}}

\newcommand{\ba}{\textbf{a}}

\newcommand{\bv}{\textbf{v}}

\newcommand{\bw}{\textbf{w}}
\newcommand{\bx}{\textbf{x}}
\newcommand{\by}{\textbf{y}}
\newcommand{\bz}{\textbf{z}}

\newcommand{\rot}{\textnormal{rot}}

\newcommand{\cP}{{\mathcal P}}

\newcommand{\mcH}{{\mathcal H}}

\newcommand{\bu}{\textnormal{$\textbf{u}$}}
\newcommand{\id}{\mathbb I}

\newcommand{\twovec}[2]{\begin{bmatrix} #1 \\ #2 \end{bmatrix}}

\newcommand{\twovecT}[2]{{\begin{bmatrix} #1 & #2 \end{bmatrix}^T}}
\newcommand{\threevecT}[3]{{\begin{bmatrix} #1 & #2 & #3 \end{bmatrix}}^T}
\newcommand{\nedelec}{N{\'e}d{\'e}lec}

\newcommand{\nl}[1]{\nabla\lambda_{#1}}
\newcommand{\lnl}[2]{\lambda_{#1}\nabla\lambda_{#2}}
\newcommand{\hlnhl}[2]{\hat\lambda_{#1}\nabla\hat\lambda_{#2}}
\newcommand{\nlxnl}[2]{\nabla\lambda_{#1}\times\nabla\lambda_{#2}}
\newcommand{\nldrnl}[2]{\nabla\lambda_{#1}\cdot \rot \nabla\lambda_{#2}}
\newcommand{\lnlxnl}[3]{\lambda_{#1}\nabla\lambda_{#2}\times\nabla\lambda_{#3}}
\newcommand{\hatlnlxnl}[3]{\hat\lambda_{#1}\nabla\hat\lambda_{#2}\times\nabla\hat\lambda_{#3}}
\newcommand{\lnldrnl}[3]{\lambda_{#1}\nabla\lambda_{#2}\cdot\rot\nabla\lambda_{#3}}
\newcommand{\lnldnlxnl}[4]{\lambda_{#1}\nabla\lambda_{#2}\cdot(\nabla\lambda_{#3}\times\nabla\lambda_{#4})}

\newtheorem{prop}[theorem]{Proposition}
\newtheorem{remark}[theorem]{Remark}
\newtheorem{cor}[theorem]{Corollary}

\begin{document}

\title{Construction of scalar and vector finite element families on polygonal and polyhedral meshes}
\author{Andrew Gillette\thanks{Department of Mathematics, University of Arizona, Tucson, AZ, USA, {\tt agillette@math.arizona.edu}}
 \and 
Alexander Rand\thanks{CD-adapco, Austin, TX, USA, {\tt alexander.rand@cd-adapco.com}}
 \and
Chandrajit Bajaj\thanks{Department of Computer Science, Institute for Computational Engineering and Sciences, University of Texas at Austin, Austin, TX, USA, {\tt bajaj@cs.utexas.edu}}
}


\maketitle

\begin{abstract}
We combine theoretical results from polytope domain meshing, generalized barycentric coordinates, and finite element exterior calculus to construct scalar- and vector-valued basis functions for conforming finite element methods on generic convex polytope meshes in dimensions 2 and 3.
Our construction recovers well-known bases for the lowest order N{\'e}d{\'e}lec, Raviart-Thomas, and Brezzi-Douglas-Marini elements on simplicial meshes and generalizes the notion of Whitney forms to non-simplicial convex polygons and polyhedra.
We show that our basis functions lie in the correct function space with regards to global continuity and that they reproduce the requisite polynomial differential forms described by finite element exterior calculus.
We present a method to count the number of basis functions required to ensure these two key properties.
\end{abstract}


\section{Introduction}

In this work, we join and expand three threads of research in the analysis of modern finite element methods: polytope domain meshing, generalized barycentric coordinates, and families of finite-dimensional solution spaces characterized by finite element exterior calculus.
It is well-known that on simplicial meshes, standard barycentric coordinates provide a local basis for the lowest-order $H^1$-conforming scalar-valued finite element spaces, commonly called the Lagrange elements.
Further, local bases for the lowest-order vector-valued Brezzi-Douglas-Marini~\cite{BDM85}, Raviart-Thomas~\cite{RT1977}, and \nedelec~\cite{BDDM87,N1980,N1986} finite element spaces on simplices can also be defined in a canonical fashion from an associated set of standard barycentric functions.
Here, we use generalized barycentric coordinates in an analogous fashion on meshes of convex polytopes, in dimensions 2 and 3, to construct local bases with the same global continuity and polynomial reproduction properties as their simplicial counterparts.

We have previously analyzed linear order, scalar-valued methods on polygonal meshes~\cite{GRB2011,RGB2011b} using four different types of generalized barycentric coordinates: Wachspress~\cite{W1975,W2011}, Sibson~\cite{F1990,S1980}, harmonic~\cite{C2008,JMRGS07,MKBWG2008}, and mean value~\cite{F2003,FHK2006,FKR2005}.
The analysis was extended by Gillette, Floater and Sukumar in the case of Wachspress coordinates to convex polytopes in any dimension~\cite{FGS2013}, based on work by Warren and colleagues~\cite{JSWD2005,W1996,WSHD2007}.
We have also shown how taking pairwise products of generalized barycentric coordinates can be used to construct quadratic order methods on polygons~\cite{RGB2011a}.
Applications of generalized barycentric coordinates to finite element methods have primarily focused on scalar-valued PDE problems~\cite{MP2008,RS2006,SM2006,ST2004,WBG07}.

\begin{table}
\centering
\begin{tabular}{c|c|c}
~n & k & functions\\ 
\hline
~2 & 0 & $\lambda_i$ \\  
  & 1 & $\lnl ij$\\
  &   & $\whit_{ij}$ \\
  &   & $\rot~\lnl ij$ \\
  &   & $\rot~\whit_{ij}$ \\
  & 2 & $\lnldrnl ijk$ \\
  &   & $\whit_{ijk}$ \\ 
\end{tabular}
\qquad
\qquad
\begin{tabular}{c|c|c}
~n & k & functions\\ 
  \hline
~3 & 0 & $\lambda_i$ \\  
  & 1 & $\lnl ij$ \\
  &   & $\whit_{ij}$ \\
  & 2 & $\lnlxnl ijk$ \\
  &   & $\whit_{ijk}$ \\
  & 3 & $\lnldnlxnl ijk\ell$   \\
  &   & $\whit_{ijk\ell}$  \\ 
\end{tabular}
\caption{For meshes of convex $n$-dimensional polytopes in $\R^n$, $n=2$ or $3$, computational basis functions for each differential form order $0\leq k\leq n$ are listed.
The notation is defined in Section~\ref{sec:bkgd}.
}
\label{tab:comp-bases}
\end{table}

Our expansion in this paper to vector-valued methods is inspired by Whitney differential forms, first defined in~\cite{W1957}.
Bossavit recognized that Whitney forms could be used to construct basis functions for computational electromagnetics~\cite{B1988a}.
The theory of finite element exterior calculus unified subsequent research in this area~\cite{AFW2006}.
In particular, Arnold, Falk and Winther showed how functions like those appearing in Table~\ref{tab:comp-bases} can be used to build spanning sets and bases for any the $\cP_r\Lambda^k$ and $\cP_r^-\Lambda^k$ spaces on simplices~\cite{AFW2009}.
The FENiCS Project~\cite{Aetal2015} has implemented these functions on simplices as part of a broadly applicable open source finite element software package.

Some prior work has explored the possibility of Whitney functions over non-simplicial elements in specific cases of rectangular grids~\cite{G2002}, square-base pyramids~\cite{GH1999}, and prisms~\cite{B2008}. 
Other authors have examined the ability of generalized Whitney functions to recover constant-valued forms in certain cases~\cite{ESW2006,KRS2011}, whereas here we show their ability to reproduce \textbf{\textit{all}} the elements of the spaces denoted $\cP_1^-\Lambda^k$ in finite element exterior calculus.
Gillette and Bajaj considered the use of generalized Whitney forms on polytope meshes defined by duality from a simplicial mesh~\cite{GB2010,GB2011}, which illustrated potential benefits to discrete exterior calculus~\cite{H2003}, computational magnetostatics, and Darcy flow modeling.
Recent work~\cite{MRS2014} has also shown generalized barycentric coordinates to be effective when used in tandem with virtual element methods~\cite{dVBCMMR2013}, which are developed in a similar fashion to traditional mimetic methods~\cite{LMS2014}.

\begin{table}[h]
\centering
\begin{tabular}{c|c|l|l}
n & k &  global continuity & polynomial reproduction\\ 
\hline
2 & 0  & $H^1(\mesh)$ & $\cP_1 \Lambda^0(\mesh)$  \\  
  & 1  & $H(\Curl,\mesh)$, by Theorem~\ref{thm:hcurl-conf} & $\cP_1 \Lambda^1(\mesh)$, by Theorem~\ref{thm:pr-lnl} \\
  &   & $H(\Curl,\mesh)$, by Theorem~\ref{thm:hcurl-conf} & $\cP_1^- \Lambda^1(\mesh)$, by Theorem~\ref{thm:pr-whit-ij}\\
  &   & $H(\Div,\mesh)$, see Remark~\ref{rmk:rot-for-cty} & $\cP_1 \Lambda^1(\mesh)$, by Corollary~\ref{cor:pr-lrnl}   \\
  &   & $H(\Div,\mesh)$, see Remark~\ref{rmk:rot-for-cty}  & $\cP_1^- \Lambda^1(\mesh)$, by Corollary~\ref{cor:pr-rotWhit} \\
  & 2  & none (piecewise linear) & $\cP_1 \Lambda^2(\mesh)$, by Theorem~\ref{thm:pr-lnldrnl}   \\
  &   & none (piecewise constant) & $\cP_1^- \Lambda^2(\mesh)$, see Remark~\ref{rmk:topdim} \\ 
  \hline
3 & 0 & $H^1(\mesh)$ & $\cP_1 \Lambda^0(\mesh)$  \\  
  & 1 & $H(\Curl,\mesh)$, by Theorem~\ref{thm:hcurl-conf} & $\cP_1 \Lambda^1(\mesh)$, by Theorem~\ref{thm:pr-lnl} \\
  &   & $H(\Curl,\mesh)$, by Theorem~\ref{thm:hcurl-conf} & $\cP_1^- \Lambda^1(\mesh)$, by Theorem~\ref{thm:pr-whit-ij}\\
  & 2 & $H(\Div,\mesh)$, by Theorem~\ref{thm:hdiv-conf} & $\cP_1 \Lambda^2(\mesh)$, by Theorem~\ref{thm:pr-lnlxnl}\\
  &   & $H(\Div,\mesh)$, by Theorem~\ref{thm:hdiv-conf} & $\cP_1^- \Lambda^2(\mesh)$, by Theorem~\ref{thm:pr-whit-ijk}\\
  & 3 & none (piecewise linear) & $\cP_1 \Lambda^3(\mesh)$, see Remark~\ref{rmk:topdim}  \\
  &   & none (piecewise constant) & $\cP_1^- \Lambda^3(\mesh)$, see Remark~\ref{rmk:topdim} \\ 
  \hline
\end{tabular}
\caption{Summary of the global continuity and polynomial reproduction properties of the spaces considered.}
\label{tab:results-summary}
\end{table}

Using the bases defined in Table~\ref{tab:comp-bases}, our main results are summarized in Table~\ref{tab:results-summary}.
On a mesh of convex $n$-dimensional polytopes in $\R^n$ with $n=2$ or $3$, we construct computational basis functions associated to the polytope elements for each differential form order $k$ as indicated.
Each function is built from generalized barycentric coordinates, denoted $\lambda_i$, and their gradients; formulae for the Whitney-like functions, denoted $\whit$, are given in Section~\ref{subset:whit-forms}.
In the vector-valued cases ($0<k<n$), we prove that the functions agree on tangential or normal components at inter-element boundaries, providing global continuity in $\hcurl$ or $\hdiv$.
The two families of polynomial differential forms that are reproduced, $\cP_r\Lambda^k$ and $\cP_r^-\Lambda^k$, were shown to recover and generalize the classical simplicial finite element spaces mentioned previously, via the theory of finite element exterior calculus~\cite{AFW2006,AFW2010}.

The outline of the paper is as follows.
In Section~\ref{sec:bkgd}, we describe relevant theory and prior work in the areas of finite element exterior calculus, generalized barycentric coordinates, and Whitney forms.
In Section~\ref{sec:global-cnty}, we show how the functions listed in Table~\ref{tab:comp-bases} can be used to build piecewise-defined functions with global continuity in $H^1$, $\hcurl$ or $\hdiv$, as indicated.
In Section~\ref{sec:poly-repro}, we show how these same functions can reproduce the requisite polynomial differential forms from $\cP_1\Lambda^k$ or $\cP_1^-\Lambda^k$, as indicated in Table~\ref{tab:comp-bases}, by exhibiting explicit linear combinations whose coefficients depend only on the location of the vertices of the mesh.
In Section~\ref{sec:polyg-fams}, we count the basis functions constructed by our approach on generic polygons and polyhedra and explain how the size of the basis could be reduced in certain cases.

\section{Background and prior work}
\label{sec:bkgd}

\subsection{Spaces from Finite Element Exterior Calculus}
\label{subsec:feec}

Finite element spaces can be broadly classified according to three parameters: $n$, the spatial dimension of the domain, $r$, the order of error decay, and $k$, the differential form order of the solution space.
The $k$ parameter can be understood in terms of the classical finite element sequence for a domain $\Omega\subset\R^n$ with $n=2$ or $3$, commonly written as
\[\xymatrix @R=.02in{
n=2: & {H^1} \ar[rr]^-{\text{grad}} &&  {\hcurl}   \ar@{<->}[rr]^-{\text{rot}} && {\hdiv} \ar[rr]^-{\text{div}}  && {L^2}\\
n=3: & {H^1} \ar[rr]^-{\text{grad}} &&  {\hcurl}   \ar[rr]^-{\text{curl}} && {\hdiv} \ar[rr]^-{\text{div}} && {L^2} 
}\]
Note that for $n=2$, given $\vec F(x,y) := \twovec{F_1(x,y)}{F_2(x,y)}$, we use the definitions:
\[\Curl\,\vec F := \frac{\p F_1}{\p y}- \frac{\p F_2}{\p x},\quad\rot\,\vec F := \begin{bmatrix} 0 & {-1} \\ 1 & 0 \end{bmatrix}\vec F\quad\text{and}\quad  \Div\vec F := \frac{\p F_1}{\p x}+ \frac{\p F_2}{\p y}.\]
\text{}\\ 
Thus, in $\R^2$, we have both $\Curl\nabla\phi =0$ and $\Div\rot\,\nabla\phi=0$ for any $\phi\in H^2$.
Put differently, $\rot$ gives an isomorphism from $\hcurl$ to $\hdiv$ in $\R^2$.
In some cases we will write $H(\Curl,\Omega)$ and $H(\Div,\Omega)$ if we wish to emphasize the domain in consideration.

In the terminology of differential topology, the applicable sequence is described more simply as the $L^2$ deRham complex of $\Omega$.
The spaces are re-cast as differential form spaces $H\Lambda^k$ and the operators as instances of the exterior derivative $d_k$, yielding
\[\xymatrix @R=.02in{
n=2: & {H\Lambda^0} \ar[rr]^-{d_0} &&  {H\Lambda^1} \ar[rr]^-{d_1} && {H\Lambda^2} \\
n=3: & {H\Lambda^0} \ar[rr]^-{d_0} &&  {H\Lambda^1}   \ar[rr]^-{d_1} && {H\Lambda^2} \ar[rr]^-{d_2} && {H\Lambda^3}
}\]

Finite element methods seek approximate solutions to a PDE in finite dimensional subspaces $\Lambda^k_h$ of the $H\Lambda^k$ spaces, where $h$ denotes the maximum diameter of a domain element associated to the subspace.
The theory of finite element exterior calculus classifies two families of suitable choices of $\Lambda^k_h$ spaces on meshes of simplices, denoted $\cP_r\Lambda^k$ and $\cP_r^-\Lambda^k$~\cite{AFW2006,AFW2010}.
The space $\cP_r\Lambda^k$ is defined as ``those differential forms which, when applied to a constant vector field, have the indicated polynomial dependence''~\cite[p. 328]{AFW2010}.
This can be interpreted informally as the set of differential $k$ forms with polynomial coefficients of total degree at most $r$.
The space $\cP_r^-\Lambda^k$ is then defined as the direct sum
\begin{equation}
\label{eq:prminus-decomp}
\cP_r^-\Lambda^k := \cP_{r-1}\Lambda^k \oplus \kappa\mcH_{r-1}\Lambda^{k+1},
\end{equation}
where $\kappa$ is the Koszul operator and $\mcH_r$ denotes homogeneous polynomials of degree $r$~\cite[p. 331]{AFW2010}.
We will use the coordinate formulation of $\kappa$, given in~\cite[p. 329]{AFW2010} as follows.
Let $\omega\in\Lambda^k$ and suppose that it can be written in local coordinates as $\omega_x=a(x)dx_{\sigma_1}\wedge\cdots\wedge dx_{\sigma_k}$.
Then $\kappa\omega$ is written as
\begin{equation}
\label{eq:def-kappa}
(\kappa\omega)_x := \sum_{i=1}^k (-1)^{i+1}a(x)x_{\sigma(i)}dx_{\sigma_1}\wedge\cdots\wedge\widehat{dx_{\sigma_i}}\wedge\cdots\wedge dx_{\sigma_k},
\end{equation}
where $\wedge$ denotes the wedge product and $\widehat{dx_{\sigma_i}}$ means that the term is omitted.
For example, let $n=3$ and write $x,y,z$ for $x_1,x_2,x_3$.
Then $dydz\in\mcH_0\Lambda^2$ and $\kappa dydz = ydz-zdy\in\mcH_1\Lambda^1$.
We summarize the relationship between the spaces $\cP_1\Lambda^k$, $\cP_1^-\Lambda^k$ and certain well-known finite element families in dimension $n=2$ or $3$ in Table~\ref{tab:fe-spaces}.

\begin{table}[ht]
\begin{center}
\begin{tabular}{c|c|c|c|l|c}
n & k & dim & space & classical description & reference\\ \hline
2 & 0 & 3 & $\cP_1 \Lambda^0(\mathcal{T})$ & Lagrange, degree $\leq 1$ &\\
  &   & 3 & $\cP_1^- \Lambda^0(\mathcal{T})$ & Lagrange, degree $\leq 1$ &\\
  & 1 & 6 & $\cP_1 \Lambda^1(\mathcal{T})$ & Brezzi-Douglas-Marini, degree $\leq 1$ & \cite{BDM85}\\
  &   & 3 & $\cP_1^- \Lambda^1(\mathcal{T})$ & Raviart-Thomas, order $0$ & \cite{RT1977}\\
  & 2 & 3 & $\cP_1 \Lambda^2(\mathcal{T})$ & discontinuous linear &\\
  &   & 1 & $\cP_1^- \Lambda^2(\mathcal{T})$ & discontinuous piecewise constant &\\ \hline
3 & 0 & 4 & $\cP_1 \Lambda^0(\mathcal{T})$ & Lagrange, degree $\leq 1$ &\\
  &   & 4 & $\cP_1^- \Lambda^0(\mathcal{T})$ & Lagrange, degree $\leq 1$ &\\
  & 1 & 12 & $\cP_1 \Lambda^1(\mathcal{T})$ & \nedelec~second kind $\hcurl$, degree $\leq 1$ & \cite{N1986,BDDM87}\\
  &   & 6 & $\cP_1^- \Lambda^1(\mathcal{T})$ & \nedelec~first kind $\hcurl$, order $0$ & \cite{N1980}\\
  & 2 & 12 & $\cP_1 \Lambda^2(\mathcal{T})$ & \nedelec~second kind $\hdiv$, degree $\leq 1$ &\cite{N1986,BDDM87}\\
  &   & 4 & $\cP_1^- \Lambda^2(\mathcal{T})$ & \nedelec~first kind $\hdiv$, order $0$ & \cite{N1980}\\
  & 3 & 4 & $\cP_1 \Lambda^3(\mathcal{T})$ & discontinuous linear & \\
  &   & 1 & $\cP_1^- \Lambda^3(\mathcal{T})$ & discontinuous piecewise constant & \\ \hline
\end{tabular}
\end{center}
\caption{Correspondence between $\cP_1\Lambda^k(\mathcal{T})$, $\cP_1^-\Lambda^k(\mathcal{T})$ and common finite element spaces associated to a simplex $\mathcal{T}$ of dimension $n$. 
Further explanation of these relationships can be found in~\cite{AFW2006,AFW2010}.
Our constructions, when reduced to simplices, recover known local bases for each of these spaces.}
\label{tab:fe-spaces}
\end{table}

A crucial property of $\cP_r\Lambda^k$ and $\cP_r^-\Lambda^k$ is that each includes in its span a sufficient number of polynomial differential $k$-forms to ensure an \textit{a priori} error estimate of order $r$ in $H\Lambda^k$ norm.
In the classical description of finite element spaces, this approximation power is immediate; any computational or `local' basis used for implementation of these spaces must, by definition, span the requisite polynomial differential forms.
The main results of this paper are proofs that generalized barycentric coordinates can be used as local bases on polygonal and polyhedral element geometries to create analogues to the lowest order $\cP_r\Lambda^k$ and $\cP_r^-\Lambda^k$ spaces with the same polynomial approximation power and global continuity properties.

In the remainder of the paper, we will frequently use standard vector proxies~\cite{AMR1988} in place of differential form notation, as indicated here:
\begin{align*}
\twovecT {u_1}{u_2} & \;\;\longleftrightarrow \;\; u_1dx_1+u_2dx_2\in\Lambda^1(\R^2), \\
\threevecT {v_1}{v_2}{v_3} & \;\;\longleftrightarrow\;\;  v_1dx_1+v_2dx_2+v_3dx_3\in\Lambda^1(\R^3),\\
\threevecT {w_1}{w_2}{w_3} & \;\;\longleftrightarrow\;\;  w_1dx_2dx_3+w_2dx_3dx_1+w_3dx_1dx_2\in\Lambda^2(\R^3).
\end{align*}

\subsection{Generalized Barycentric Coordinates}
\label{subsec:gbcs}

Let $\elt$ be a convex $n$-dimensional polytope in $\R^n$ with vertex set $\{\bv_i\}$, written as column vectors.
A set of non-negative functions $\{\lambda_i\}:\elt\raw\R$ are called \textbf{generalized barycentric coordinates} on $\elt$ if for any linear function $L:\elt\raw\R$, we can write
\begin{equation}
\label{eq:lin-comp}
L = \sum_i L(\bv_i)\lambda_i,
\end{equation}
We will use the notation $\id$ to denote the $n\times n$ identity matrix and $\bx$ to denote the vector ${\begin{bmatrix} x_1 & x_2 & \cdots & x_n\end{bmatrix}}^T$ where $x_i$ is the $i$th coordinate in $\R^n$.
We have the following useful identities:
\begin{align}
\sum_i\lambda_i(\bx) &= 1 \label{eq:pof1} \\
\sum_i\bv_i\lambda_i(\bx) &= \bx \label{eq:pofx}\\
\sum_i\nabla\lambda_i (\bx)& = 0 \label{eq:gradsum0} \\
\sum_i\bv_i\nabla\lambda_i^T(\bx) & = \id \label{eq:vgradsumI}
\end{align}
Equations (\ref{eq:pof1}) and (\ref{eq:pofx}) follow immediately from (\ref{eq:lin-comp}) while (\ref{eq:gradsum0}) and (\ref{eq:vgradsumI}) follow by taking the gradient of equations (\ref{eq:pof1}) and (\ref{eq:pofx}), respectively.
If $\bx$ is constrained to an $n-1$ dimensional facet of $\elt$ and the index set of the summations are limited to those vertices that define $\elt$, then (\ref{eq:pof1})-(\ref{eq:vgradsumI}) still hold; in particular, this implies that generalized barycentric coordinates on a polyhedron restrict to generalized barycentric coordinates on each of its polygonal faces.

As mentioned in the introduction, there are many approaches to defining generalized barycentric coordinates.
In regards to applications in finite element methods, the Wachspress coordinates~\cite{W1975,W2011} are commonly used as they are rational functions in both 2D and 3D with explicit formulae; code for their implementation in MATLAB is given in the appendix of~\cite{FGS2013}.
Other practical choices of generalized barycentric coordinates for finite elements include mean value~\cite{F2003}, maximum entropy~\cite{HS2008,Su04}, and moving least squares~\cite{MS10}.
The results of this work do not rely on any properties of the coordinates other than their non-negativity and linear reproduction property (\ref{eq:lin-comp}).

\subsection{Whitney forms}
\label{subset:whit-forms}

Let $\elt$ be a convex $n$-dimensional polytope in $\R^n$ with vertex set \textnormal{$\{\bv_i\}$} and an associated set of generalized barycentric coordinates $\{\lambda_i\}$.
Define associated sets of index pairs and triples by
\begin{align}
E & := \{(i,j)   ~:~ \bv_i,\bv_j\in\elt\} \label{eq:Em-def},\\
T & := \{(i,j,k) ~:~ \bv_i,\bv_j,\bv_k\in\elt\} \label{eq:Tm-def}.
\end{align}
If $\elt$ is a \textit{simplex}, the elements of the set
\[\left\{\lnl ij - \lnl ji~:~(i,j)\in E\right\}\]
are called Whitney 1-forms and are part of a more general construction~\cite{W1957}, which we now present.
Again, if $\elt$ is a $\textit{simplex}$, the Whitney $k$-forms are elements of the set
\begin{equation}
\label{eq:whit-def}
\left\{k!\sum_{i=0}^k(-1)^i\;  \lambda_{j_i} \;d\lambda_{j_0}\wedge\ldots\wedge\widehat{d\lambda_{j_i}}\wedge\ldots\wedge d\lambda_{j_k}\right\},
\end{equation}
where $j_0,\ldots,j_k$ are indices of vertices of $\elt$.
As before, $\wedge$ denotes the wedge product and $\widehat{dx_{\sigma_i}}$ means that the term is omitted.
Up to sign, this yields a set of $n+1\choose k+1$ distinct functions and provides a local basis for $\cP_1^-\Lambda^k$~\cite{AFW2009}.

We now generalize these definitions to the case where $\elt$ is non necessarily a simplex.
For any $(i,j)\in E$, define a generalized Whitney 1-form on $\elt$ by 
\begin{align}
\label{eq:whit-edge-def}
\whit_{ij} & := \lambda_i\nabla\lambda_j-\lambda_j\nabla\lambda_i.
\end{align}
If $n=3$, then for any $(i,j,k)\in T$, define a generalized Whitney 2-form on $\elt$ by
\begin{align}
\label{eq:whit-tri-def}
\whit_{ijk} & := (\whitC ijk) + (\whitC jki) + (\whitC kij).
\end{align}
Note that $\whit_{ii}=0$ and if $i$, $j$, and $k$ are not distinct then $\whit_{ijk}=0$.

Whitney forms have natural interpretations as vector fields when $k=1$ or $n-1$.
Interpolation of vector fields requires less data regularity than the canonical scalar interpolation theory using nodal values.
Averaged interpolation developed for scalar spaces~\cite{Cl75,SZ90} has been extended to families of spaces from finite element exterior calculus~\cite{CW08}. 
Recent results on polygons and polyhedra can be extended to less regular data with average interpolation following the framework in~\cite{Ra12}, based on affine invariance of the coordinates.

\section{Global Continuity Results}
\label{sec:global-cnty}

We first present results about the global continuity properties of vector-valued functions defined in terms of generalized barycentric coordinates and their gradients over a mesh of $n$-dimensional polytopes in $\R^n$ with $n=2$ or $3$.
By `mesh' we mean a cellular complex in which each cell is a polygon (for $n=2$) or polyhedron (for $n=3$); for more on cellular complexes see e.g.~\cite{C2008}.
Voronoi meshes are examples of cellular complexes since they are composed of $n$-dimensional polytopes that meet along their $n-1$ dimensional facets.
We say that a function is defined `piecewise with respect to a mesh' when the definition of the function on the interior of a mesh element depends only on geometrical properties of the element (as opposed to depending on adjacent elements, for instance).
We begin with a general result about global continuity in such a setting.\\

\begin{prop}
\label{prop:trace2conf}
Fix a mesh $\mesh$ of $n$-dimensional polytopes in $\R^n$ with $n=2$ or $3$.
Let $\bu$ be a vector field defined piecewise with respect to $\mesh$.
Let $\face$ be a face of codimension 1 with $\bu_1$, $\bu_2$ denoting the values of $\bu$ on $\face$ as defined by the two $n$-dimensional mesh elements sharing $\face$. 
Write $\bu_i = T_\face(\bu_i) + N_\face(\bu_i)$ where $T_\face(\bu_i)$ and $N_\face(\bu_i)$ are the vector projections of $\bu_i$ onto $\face$ and its outward normal, respectively. 
\renewcommand{\labelenumi}{(\roman{enumi}.)}
\begin{enumerate}
\item If $T_\face(\bu_1)=T_\face(\bu_2)$ for all $\face\in \mesh$ then $\bu\in H(\Curl,\mesh)$.
\item If $N_\face(\bu_1)=N_\face(\bu_2)$ for all $\face\in \mesh$ then $\bu\in H(\Div,\mesh)$.
\end{enumerate}
\end{prop}
\text{}\\
The results of Proposition~\ref{prop:trace2conf} are well-known in the finite element community; see e.g.~Ern and Guermond~\cite[Section 1.4]{EG04}.\\

\begin{prop}
\label{prop:bctrcfaceprop}
Let $\elt$ be a convex $n$-dimensional polytope in $\R^n$ with vertex set \textnormal{$\{\bv_i\}_{i\in I}$} and an associated set of generalized barycentric coordinates $\{\lambda_i\}_{i\in I}$.
Let $\face$ be a face of $\elt$ of codimension 1 whose vertices are indexed by $J\subsetneq I$.
If $k\not\in J$ then $\lambda_k\equiv 0$ on $\face$ and $\nabla\lambda_k$ is normal to $\face$ on $\face$, pointing inward.
\end{prop}
\begin{proof}
Fix a point $\bx_0\in\elt$.
Observe that $\sum_{i\in I}\bv_i\lambda_i(\bx_0)$ is a point in $\elt$ lying in the interior of the convex hull of those $\bv_i$ for which $\lambda_i(\bx_0)>0$, since the $\lambda_i$ are non-negative by definition.
By (\ref{eq:pofx}), this summation is equal to $\bx_0$.
Hence, if $\bx_0\in \face$, then $\lambda_k\equiv 0$ on $\face$ unless $k\in J$, proving the first claim.
The same argument implies that for any $k\not\in $J, $\face$ is part of the zero level set of $\lambda_k$.
Hence, for $k\not\in J$, $\nabla\lambda_k$ is orthogonal to $\face$ on $\face$.  
In that case, $\nabla\lambda_k$ points inward since $\lambda_k$ has support inside $\elt$ but not on the other side of $\face$.
\end{proof}

\begin{figure}
\begin{center}
\includegraphics[height=.3\textwidth]{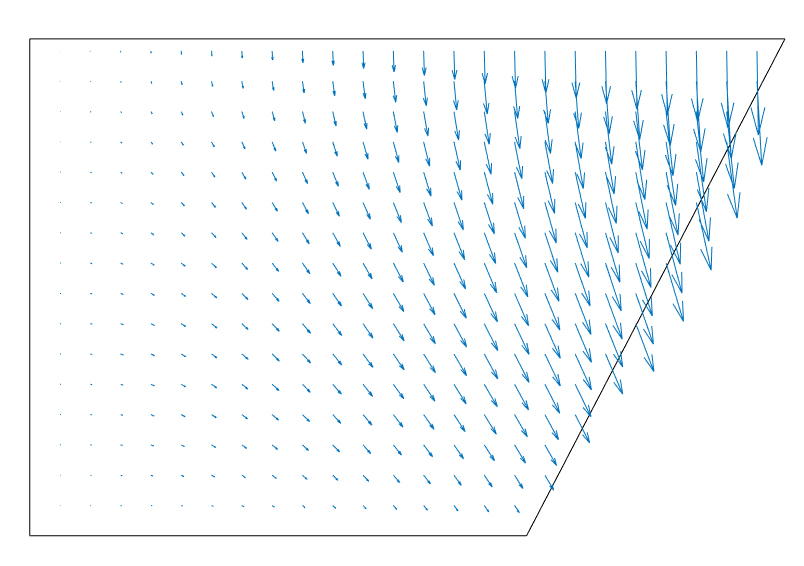} ~\includegraphics[height=.3\textwidth]{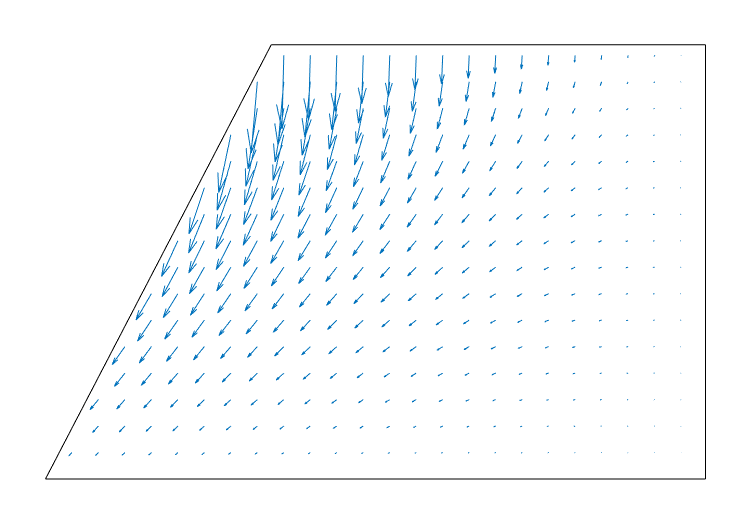}
\end{center}
\caption{The $\hcurl$ conformity condition of Proposition~\ref{prop:trace2conf} is satisfied automatically by the $\lambda_i\nabla\lambda_j$ functions, as shown in the example above.
When the elements are brought together, the vector fields will agree on the projection to the shared edge at any point along the shared edge.
Here, $i$ and $j$ are the indices for the vertices at the top and bottom, respectively, of the shared edge.  For this example, we used the Wachspress functions to compute the vector functions on each element and MATLAB to visualize the result.}
\label{fig:conf}
\end{figure}

We now show that generalized barycentric coordinates and their gradients defined over individual elements in a mesh of polytopes naturally stitch together to build conforming finite elements with global continuity of the expected kind.
Figure~\ref{fig:conf} presents an example of two vector functions agreeing on their tangential projections along a shared edge.
To be clear about the context, we introduce notation for generalized barycentric hat functions, defined piecewise over a mesh of polytopes $\{\elt\}$ by
\[\hat\lambda_i(\bx)= \begin{cases} \lambda_i(\bx)~\text{as defined on $\elt$} & \text{if $\bx\in\elt$ and $\bv_i\in \elt$;} \\ 0 & \text{if $\bx\in\elt$ but $\bv_i \not\in \elt$.} \end{cases}\]
Note that generalized barycentric coordinates $\lambda_i$ are usually indexed locally on a particular polytope while the $\hat\lambda_i$ require a global indexing of the vertices to consistently identify matching functions across element boundaries. 
Further, $\hat\lambda_i$ is well-defined at vertices and edges of the mesh as any choice of generalized barycentric coordinates on a particular element will give the same value at such points.
If $\bx$ belongs to the interior of shared faces between polyhedra in $\R^3$ (or higher order analogues), $\hat\lambda_i(\bx)$ is well-defined so long as the same \textit{kind} of generalized barycentric coordinates are chosen on each of the incident polyhedra (e.g.\ Wachspress or mean value).

Our first result about global continuity concerns functions of the form $\hlnhl ij$, where $i$ and $j$ are indices of vertices belonging to at least one fixed mesh element $\elt$.
Note that the vertices $\bv_i$ and $\bv_j$ need not define an edge of $\elt$.\\

\begin{theorem}
\label{thm:hcurl-conf}
Fix a mesh $\mesh$ of $n$-dimensional polytopes $\{\elt\}$ in $\R^n$ with $n=2$ or $3$ and assign some ordering $\textnormal{$\bv_1,\ldots,\bv_p$}$ to all the vertices in the mesh.
Fix an associated set of generalized barycentric coordinate hat functions $\textnormal{$\hat\lambda_1,\ldots,\hat\lambda_p$}$.
Let
\[\bu\in\textnormal{span}~\left\{\hlnhl ij~:~\textnormal{$\exists~\elt\in\mesh$ such that $\bv_i,\bv_j\in\elt$}\right\}.\]
Then $\bu\in H(\Curl,\mesh)$.
\end{theorem}

\begin{proof}
Following the notation of Proposition~\ref{prop:trace2conf}, it suffices to show that  $T_\face(\bu_1)=T_\face(\bu_2)$ for an arbitrary face $\face\in \mesh$ of codimension 1.
Consider an arbitrary term $c_{ij}\hlnhl ij$ in the linear combination defining $\bu$.
Observe that if $\bv_i\not\in \face$, then by Proposition~\ref{prop:bctrcfaceprop}, $\hat \lambda_i\equiv 0$ on $\face$ and hence $\bu\equiv 0$ on $\face$.
Further, if $\bv_j\not\in \face$, then $\nabla\hat \lambda_j$ is orthogonal to $\face$.
Therefore, without loss of generality, we can reduce to the case where $\bv_i,\bv_j\in \face$.
Since $\hat \lambda_i$ and $\hat \lambda_j$ are both $C^0$ on $\mesh$, their well-defined values on $\face$ suffice to determine the projection of $\hlnhl ij$ to $\face$.
Since the choice of pair $ij$ was arbitrary, we have $T_\face(\bu_1)=T_\face(\bu_2)$, completing the proof.
\end{proof}

\text{}
\begin{remark}
\label{rmk:rot-for-cty}
{\em
When $n=2$, we may replace $\hlnhl ij$ in the statement Theorem~\ref{thm:hcurl-conf} by $\rot~\hlnhl ij$ and conclude that $\bu\in H(\Div,\mesh)$.
This is immediate since $\rot$ gives an isomorphism between $\hcurl$ and $\hdiv$ in $\R^2$, as discussed in Section~\ref{subsec:feec}.
When $n=3$, we construct functions in $H(\Div,\mesh)$ using triples of indices associated to vertices of mesh elements, according to the next result.
}
\end{remark}

\text{}
\begin{theorem}
\label{thm:hdiv-conf}
Fix a mesh $\mesh$ of polyhedra $\{\elt\}$ in $\R^3$ and assign some ordering $\textnormal{$\bv_1,\ldots,\bv_p$}$ to all the vertices in the mesh.
Fix an associated set of generalized barycentric coordinate hat functions $\textnormal{$\hat\lambda_1,\ldots,\hat\lambda_p$}$.
Let
\[\bu\in\textnormal{span}~\left\{\hatlnlxnl ijk~:~\textnormal{$\exists~\elt\in\mesh$ such that  $\bv_i,\bv_j,\bv_k\in\elt$}\right\}.\]
Then $\bu\in H(\Div,\mesh)$.
\end{theorem}

\begin{proof}
Again following the notation of Proposition~\ref{prop:trace2conf}, it suffices to show that  $N_\face(\bu_1)=N_\face(\bu_2)$ for an arbitrary face $\face\in \mesh$ of codimension one whose vertices are indexed by $J$.
We will use the shorthand notation
\[\xi_{ijk} := \hatlnlxnl ijk. \]
Consider an arbitrary term $c_{ijk}\xi_{ijk}$ in the linear combination defining $\bu$.
We will first show that $\xi_{ijk}$ has a non-zero normal component on $\face$ only if $i,j,k\in J$.
If $i\not\in J$ then $\hat\lambda_i\equiv0$ on $\face$ by Proposition~\ref{prop:bctrcfaceprop}, making $\xi_{ijk}\equiv 0$ on $\face$, as well.  
If $i\in J$ but $j,k\not\in J$, then $\nabla\hat\lambda_j$ and $\nabla\hat\lambda_k$ are both normal to $\face$ on $\face$ by Proposition~\ref{prop:bctrcfaceprop}.  
Hence, their cross product is zero and again $\xi_{ijk}\equiv 0$ on $F$.  
If $i,j\in J$ but $k\not\in J$ then again $\nabla\hat\lambda_k\perp \face$ on $\face$.  
Since $\nabla\hat\lambda_j\times\nabla\hat\lambda_k\perp \nabla\hat\lambda_k$, we conclude that $\xi_{ijk}$ has no normal component on $\face$.  
The same argument holds for the case $i,k\in J$, $j\not\in J$.  
The only remaining case is $i,j,k\in J$, proving the claim.

Thus, without loss of generality, we assume that $i,j,k\in J$.
Since $\hat\lambda_j$ and $\hat\lambda_k$ are both $C^0$ on $\mesh$, their well-defined values on $\face$ suffice to determine the projection of $\nabla\hat\lambda_j$ and $\nabla\hat\lambda_k$ to $\face$, which then uniquely defines the normal component of $\nabla\hat\lambda_j\times\nabla\hat\lambda_k$ on $\face$.
Since $\hat\lambda_i$ is also $C^0$ on $\mesh$, and the choice of $i,j,k$ was arbitrary, we have $N_\face(\bu_1)=N_\face(\bu_2)$, completing the proof.
\end{proof}
\text{}\\

\section{Polynomial Reproduction Results}
\label{sec:poly-repro}

We now show how generalized barycentric coordinate functions $\lambda_i$ and their gradients can reproduce all the polynomial differential forms in $\cP_1\Lambda^k$ and $\cP_1^-\Lambda^k$ for $0\leq k\leq n$ with $n=2$ or $3$. 
The results for the functions $\lnl ij$ and $\whit_{ij}$ extend immediately to any value of $n\geq 2$ since those functions do not use any dimension-specific operators like $\times$ or $\rot$. \\

\begin{theorem}
\label{thm:pr-lnl}
Fix $n\geq 2$.
Let $\elt$ be a convex $n$-dimensional polytope in $\R^n$ with vertex set \textnormal{$\{\bv_i\}$}.
Given any set of generalized barycentric coordinates $\{\lambda_i\}$ associated to $\elt$,
\begin{equation}
\label{eq:lnl-id}
\textnormal{$
\sum_{i,j}\lnl ij(\bv_j-\bv_i)^T= \id,
$}
\end{equation}
where $\id$ is the $n\times n$ identity matrix.
Further, for any $n\times n$ matrix $\A$,
\begin{equation}
\label{eq:lnl-Ax}
\textnormal{$
\sum_{i,j}(\A\bv_i \cdot\bv_j)(\lnl ij) =\A\bx.
$}
\end{equation}
Thus, $\ds\textnormal{span}\left\{\lnl ij\;:\;\textnormal{$\bv_i,\bv_j\in\elt$}\right\}\supseteq \cP_1\Lambda^1(\elt).$
\end{theorem}
\begin{proof}
From (\ref{eq:pof1}) - (\ref{eq:vgradsumI}), we see that
\begin{align*}
\sum_{i,j}\lnl ij(\bv_j-\bv_i)^T 
& = \left(\sum_{i}\lambda_i\right)\left(\sum_j\nabla\lambda_j\bv_j^T\right) - \left(\sum_{j}\nabla\lambda_j\right)\left(\sum_i\lambda_i\bv_i^T\right) \\
& = 1(\id^T) - 0(\bx^T) = \id,
\end{align*}
establishing (\ref{eq:lnl-id}).
Similarly for (\ref{eq:lnl-Ax}), a bit of algebra yields
\begin{align*}
\sum_{i,j}(\A\bv_i \cdot\bv_j)(\lnl ij) 
& = \sum_{i,j}(\lnl ij) \bv_j^T \A\bv_i
 = \sum_{i,j} \nabla\lambda_j \bv_j^T \A\bv_i  \lambda_i\\
& = \left(\sum_{j} \nabla\lambda_j \bv_j^T\right)\A\left(\sum_i \bv_i \lambda_i\right)
 = \id^T\A\bx
 = \A\bx
\end{align*}
We have shown that any vector of linear polynomials can be written as a linear combination of $\lnl ij$ functions, hence the span of these functions contains the vector proxies for all elements of $\cP_1\Lambda^1(\elt)$.
\end{proof}
\text{}\\

\begin{cor}
\label{cor:pr-lrnl}
Let $\elt$ be a convex polygon in $\R^2$ with vertex set \textnormal{$\{\bv_i\}$}.
Given any set of generalized barycentric coordinates $\{\lambda_i\}$ associated to $\elt$,
\begin{equation}
\label{eq:lrnl-id}
\textnormal{$
\sum_{i,j}\rot\,\lnl ij(\rot(\bv_j-\bv_i))^T= \id,
$}
\end{equation}
where $\id$ is the $2\times 2$ identity matrix.
Further, for any $2\times 2$ matrix $\A$,
\begin{equation}
\label{eq:lrnl-Ax}
\textnormal{$
\sum_{i,j}(-\rot\,\A\,\bv_i \cdot\bv_j)(\rot\, \lnl ij) =\A\bx.
$}
\end{equation}
Thus, $\ds\textnormal{span}\left\{\rot\lnl ij\;:\;\textnormal{$\bv_i,\bv_j\in\elt$}\right\}\supseteq \cP_1\Lambda^1(\elt).$
\end{cor}

\begin{proof}
For (\ref{eq:lrnl-id}), observe that for any $\bw,\by\in\R^2$, $\ds\bw\by^T = \begin{bmatrix} a & b \\ c & d\end{bmatrix}$ implies $\ds(\rot~\bw)(\rot~\by)^T = \begin{bmatrix} d & -c \\ -b & a \end{bmatrix}$.
Hence, the result follows immediately from (\ref{eq:lnl-id}).
For (\ref{eq:lrnl-Ax}), note $\rot^{-1}=-\rot$ and define $\B:=-\rot\,\A$.
Using $\B$ as the matrix in (\ref{eq:lnl-Ax}), we have
\[\sum_{i,j}(\B\bv_i \cdot\bv_j)(\lnl ij) =\B\bx\]
Applying $\rot$ to both sides of the above yields the result.
\end{proof}
\text{}\\

\begin{theorem}
\label{thm:pr-lnlxnl}
Let $\elt$ be a convex polyhedron in $\R^3$ with vertex set \textnormal{$\{\bv_i\}$}.
Given any set of generalized barycentric coordinates $\{\lambda_i\}$ associated to $\elt$,
\begin{equation}
\label{eq:lnlxnl-id}
\textnormal{$
\frac 12\sum_{i,j,k}\lnlxnl ijk\left((\bv_j-\bv_i)\times(\bv_k-\bv_i)\right)^T= \id ,
$}
\end{equation}
where $\id$ is the $n\times n$ identity matrix.
Further, for any $n\times n$ matrix $\A$,
\begin{equation}
\label{eq:lnlxnl-Ax}
\textnormal{$
\frac 12\sum_{i,j,k}(\A\bv_i \cdot(\bv_j\times\bv_k))(\lnlxnl ijk) =\A\bx.
$}
\end{equation}
Thus, $\ds\textnormal{span}\left\{\lnlxnl ijk\;:\;\textnormal{$\bv_i,\bv_j,\bv_k\in\elt$}\right\}\supseteq \cP_1\Lambda^2(\elt).$
\end{theorem}
\begin{proof}
We start with (\ref{eq:lnlxnl-id}).
First, observe that
\[(\bv_j-\bv_i)\times(\bv_k-\bv_i) = \bv_i\times\bv_j +\bv_j\times\bv_k + \bv_k\times\bv_i.\]
By (\ref{eq:gradsum0}), we have that
\begin{align*}
\sum_{i,j,k}\lnlxnl ijk\left(\bv_i\times\bv_j\right)^T  & = \sum_{i,j} \lambda_i \left(\nabla \lambda_j \times \left( \sum_k \nabla \lambda_k \right)\right) \left(\bv_i\times\bv_j\right)^T 
=0.
\end{align*}
A similar argument shows that replacing $\bv_i\times\bv_j$ with $\bv_k\times\bv_i$ also yields the zero matrix.
Hence,
\begin{align*}
\sum_{i,j,k}\lnlxnl ijk\left((\bv_j-\bv_i)\times(\bv_k-\bv_i)\right)^T &= \sum_{i,j,k}\lnlxnl ijk\left(\bv_j\times\bv_k\right)^T
\end{align*}
\begin{align*}
= \sum_i \lambda_i \sum_{j,k}\left(\nabla \lambda_j \times \nabla \lambda_k\right) \left(\bv_j\times\bv_k\right)^T
& =\sum_{j,k}\left(\nabla \lambda_j \times \nabla \lambda_k\right) \left(\bv_j\times\bv_k\right)^T.
\end{align*}
To simplify this further, we use the Kronecker delta symbol $\delta_{i_1i_2}$ and the 3D Levi-Civita symbol $\veps_{i_1 i_2 i_3}$.
It suffices to show that the entry in row $r$, column $c$ of the matrix $\sum_{j,k}(\nlxnl jk)\left(\bv_j\times \bv_k\right)^T$ is $2\delta_{rc}$.
We see that
\begin{align*}
\left[\sum_{j,k}(\nlxnl jk)\left(\bv_j\times \bv_k\right)^T\right]_{rc}  
& =  \sum_{j,k}\veps_{r\ell m}(\nl j)_\ell(\nl k)_m \veps_{cpq} (\bv_j)_p(\bv_k)_q \\
& = \veps_{r\ell m}\veps_{cpq} \sum_{j}(\bv_j)_p(\nl j)_\ell \sum_{k} (\bv_k)_q(\nl k)_m   \\
& = \veps_{r\ell m}\veps_{cpq} \delta_{\ell p}\delta_{mq}.
\end{align*}
The last step in the above chain of equalities follows from (\ref{eq:vgradsumI}).
Observe that $\veps_{r\ell m}\veps_{cpq} \delta_{\ell p}\delta_{mq}= \veps_{r\ell m}\veps_{c\ell m} = 2\delta_{r c}$, as desired.
For (\ref{eq:lnlxnl-Ax}), observe that
\begin{align*}
\sum_{i,j,k}(\A\bv_i \cdot(\bv_j\times\bv_k))(\lnlxnl ijk) 
& = \left(\sum_{i}\A\bv_i\lambda_i\right)\cdot \sum_{j,k} (\bv_j\times\bv_k)(\nlxnl jk) \\
& = \sum_{j,k} (\nlxnl jk)(\bv_j\times\bv_k)^T \left(\A \sum_{i}\bv_i\lambda_i\right) \\
& = 2\;\id\;\A\bx = 2\A\bx.
\end{align*}
Note that we used the proof of (\ref{eq:lnlxnl-id}) to rewrite the sum over $j,k$ as $2\id$.
We have shown that any vector of linear polynomials can be written as a linear combination of $\lnlxnl ijk$ functions, hence the span of these functions contains the vector proxies for all elements of $\cP_1\Lambda^2(\elt)$.
\end{proof}
\text{}\\

\begin{theorem}
\label{thm:pr-lnldrnl}
Let $\elt$ be a convex polygon in $\R^2$ with vertex set \textnormal{$\{\bv_i\}$}.
Given any set of generalized barycentric coordinates $\{\lambda_i\}$ associated to $\elt$,
\begin{equation}
\label{eq:lnldrnl-one}
\textnormal{$
\frac 12\sum_{i,j,k}\lnldrnl ijk\left((\bv_j-\bv_i)\cdot \rot(\bv_k-\bv_i)\right)= 1. 
$}
\end{equation}
Further, for any vector $\ba \in \R^2$,
\begin{equation}
\label{eq:lnldrnl-Ax}
\textnormal{$
\frac 12\sum_{i,j,k}(\ba^T\bv_i (\bv_j\cdot \rot\bv_k))(\lnldrnl ijk) =\ba^T\bx.
$}
\end{equation}
Thus, $\ds\textnormal{span}\left\{\lnldrnl ijk\;:\;\textnormal{$\bv_i,\bv_j,\bv_k\in\elt$}\right\}\supseteq \cP_1\Lambda^2(\elt).$
\end{theorem}

\begin{proof}
The proof is essentially identical to that of Theorem~\ref{thm:pr-lnlxnl}. First,
\[(\bv_j-\bv_i) \cdot \rot(\bv_k-\bv_i) = \bv_i\cdot \rot\bv_j +\bv_j\cdot \rot\bv_k + \bv_k\cdot \rot\bv_i,\]
and by (\ref{eq:gradsum0}),
\begin{align*}
\sum_{i,j,k}\lnldrnl ijk & \left(\bv_i\cdot \rot\bv_j\right) = \sum_{i,j} \lambda_i \left(\nabla \lambda_j \cdot \rot \left( \sum_k \nabla \lambda_k \right)\right) \left(\bv_i\cdot \rot \bv_j\right)
=0.
\end{align*}
A similar argument shows that replacing $\bv_i\cdot \rot \bv_j$ with $\bv_k\cdot \rot\bv_i$ also yields zero.
Hence as before,
\begin{align*}
\sum_{i,j,k}\lnldrnl ijk & \left((\bv_j-\bv_i)\cdot \rot(\bv_k-\bv_i)\right)^T  =\sum_{j,k}\left(\nabla \lambda_j \cdot \rot \nabla \lambda_k\right) \left(\bv_j\cdot \rot\bv_k\right)^T.
\end{align*}
Finally, the same argument holds using the 2D Levi-Civita symbol:
\begin{align*}
\sum_{j,k}(\nldrnl jk)\left(\bv_j\cdot\rot \bv_k\right)  
& =  \sum_{j,k}\veps_{\ell m}(\nl j)_\ell(\nl k)_m \veps_{p q} (\bv_j)_p(\bv_k)_q \\
& = \veps_{\ell m}\veps_{pq} \sum_{j}(\bv_j)_p(\nl j)_\ell \sum_{k} (\bv_k)_q(\nl k)_m   \\
& = \veps_{\ell m}\veps_{pq} \delta_{\ell p}\delta_{mq}
  = \veps_{\ell m}\veps_{\ell m} = 2,
\end{align*}
establishing (\ref{eq:lnldrnl-one}). For (\ref{eq:lnldrnl-Ax}), observe that
\begin{align*}
\sum_{i,j,k}(\ba^T\bv_i & (\bv_j\cdot \rot \bv_k)) (\lnldrnl ijk) \\
& = \left(\sum_{i}\ba^T\bv_i\lambda_i\right) \sum_{j,k} (\bv_j\cdot \rot \bv_k)(\nldrnl jk) \\
& = \sum_{j,k} (\nldrnl jk)(\bv_j\cdot \rot \bv_k)^T \left(\ba^T \sum_{i}\bv_i\lambda_i\right)
 = 2\ba^T\bx.
\end{align*}
\end{proof}

\begin{remark}
{\em
The proof of Theorem~\ref{thm:pr-lnldrnl} can also be obtained by augmenting the 2D vectors and matrices with zeros to make 3D vectors and matrices and recognizing (\ref{eq:lnldrnl-one}) as the element equality in the third row and third column of (\ref{eq:lnlxnl-id}).
}
\end{remark}
\text{}\\

We also have polynomial reproduction results using the Whitney-like basis functions (\ref{eq:whit-edge-def}) and (\ref{eq:whit-tri-def}).
Recall that $\mcH_r$ denotes homogeneous polynomials of degree $r$ and let $\M_{n\times n}$ denote $n\times n$ matrices. 
We have the following theorems.\\

\begin{theorem}
\label{thm:pr-whit-ij}
Fix $n\geq 2$.
Let $\elt$ be a convex $n$-dimensional polytope in $\R^n$ with vertex set \textnormal{$\{\bv_i\}$} and an associated set of generalized barycentric coordinates $\{\lambda_i\}$.
Then
\begin{equation}
\label{eq:whit1-id}
\textnormal{$
\sum_{i<j}\whit_{ij}(\bv_j-\bv_i)^T=\id.
$}
\end{equation}
Further, define a map $\Phi:\mcH_1\Lambda^1(\R^n)\raw \M_{n\times n}$ by
\[\textnormal{$
\sum_{i=1}^n\left(\sum_{j=1}^n a_{ij}x_j\right) dx_i \longmapsto \left [ \sign(a_{ij}) \right ] .
$}\]
Then for all $\omega\in \mcH_0\Lambda^2(\R^n)$, 
\begin{equation}
\label{eq:whit1-x}
\textnormal{$
\sum_{i<j}\left(\Phi(\kappa\omega)\bv_i)\cdot \bv_j\right) \whit_{ij} = (\Phi(\kappa\omega))\bx.
$}
\end{equation}
Thus, $\ds\textnormal{span}\left\{\whit_{ij}\;:\;\textnormal{$\bv_i,\bv_j\in\elt$}\right\}\supseteq \cP_1^-\Lambda^1(\elt).$
\end{theorem}
\begin{proof}
For (\ref{eq:whit1-id}), we reorganize the summation and apply (\ref{eq:lnl-id}) to see that
\begin{align*}
\sum_{i<j}\whit_{ij}(\bv_j-\bv_i)^T  
 & = \sum_{i<j}\lnl ij(\bv_j-\bv_i)^T - \sum_{i<j}\lnl ji(\bv_j-\bv_i)^T \\
 & = \sum_{i<j}\lnl ij(\bv_j-\bv_i)^T + \sum_{j<i}\lnl ij(\bv_j-\bv_i)^T \\
 & =\sum_{i,j}\lnl ij(\bv_j-\bv_i)^T = \id.
\end{align*}
For (\ref{eq:whit1-x}), fix $\omega\in \mcH_0\Lambda^2(\R^n)$ and express it as
\[\omega = \sum_{i<j}a_{ij}dx_idx_j,\]
for some coefficients $a_{ij}\in\R$.
Then 
\[\kappa\omega = \sum_{i<j}a_{ij}(x_idx_j-x_jdx_i).\]
The entries of the matrix $\Phi(\kappa\omega)$ are thus given by
\begin{equation}
\label{eq:phi-entries}
\left[\Phi(\kappa\omega)\right]_{ij} = 
\begin{cases} 
\sign(a_{ij}) & \mbox{if } i<j, \\
-\sign(a_{ij}) & \mbox{if } i>j, \\
0 & \mbox{if } i=j. \\
\end{cases}
\end{equation}
From (\ref{eq:lnl-Ax}), we have that
\[\sum_{i,j}\left(\Phi(\kappa\omega)\bv_i)\cdot \bv_j\right) \lnl ij = (\Phi(\kappa\omega))\bx,\qquad \forall \omega\in \mcH_0\Lambda^2(\R^n)\]
Since $\Phi(\kappa\omega)$ is anti-symmetric by (\ref{eq:phi-entries}), we have that
\begin{align*}
\sum_{i,j}\left(\Phi(\kappa\omega)\bv_i)\cdot \bv_j\right) & \lnl ij  \\
& = \sum_{i<j}\left(\Phi(\kappa\omega)\bv_i)\cdot \bv_j\right) \lnl ij + \sum_{j<i}\left(\Phi(\kappa\omega)\bv_i)\cdot \bv_j\right) \lnl ij  \\
& = \sum_{i<j}\left(\Phi(\kappa\omega)\bv_i)\cdot \bv_j\right) \whit_{ij}.
\end{align*}
We have shown that any vector proxy of an element of $\cP_0\Lambda^1(\elt)$ or $\kappa\mcH_{0}\Lambda^{2}(\elt)$ can be written as a linear combination of $\whit_{ij}$ functions.
By (\ref{eq:prminus-decomp}), we conclude that the span of the $\whit_{ij}$ functions  contains the vector proxies for all elements of $\cP_1^-\Lambda^1(\elt)$.
\end{proof}
\text{}\\

\begin{cor}
\label{cor:pr-rotWhit}
Let $\elt$ be a convex polygon in $\R^2$ with vertex set \textnormal{$\{\bv_i\}$}.
Given any set of generalized barycentric coordinates $\{\lambda_i\}$ associated to $\elt$,
\begin{equation}
\label{eq:rotWhit-id}
\textnormal{$
\sum_{i<j}\rot~\whit_{ij}~\rot(\bv_j-\bv_i)^T=\id,
$}
\end{equation}
where $\id$ is the $2\times 2$ identity matrix.
Further, 
\begin{equation}
\label{eq:rotWhit-x}
\textnormal{$
\sum_{i<j}\left((\rot~\bv_i)\cdot \bv_j\right) \rot~\whit_{ij} = \bx.
$}
\end{equation}
Thus, $\ds\textnormal{span}\left\{\rot\;\whit_{ij}\;:\;\textnormal{$\bv_i,\bv_j\in\elt$}\right\}\supseteq \cP_1^-\Lambda^1(\elt).$
\end{cor}

\begin{proof}
By the same argument as the proof of (\ref{eq:lrnl-id}) in Corollary~\ref{cor:pr-lrnl}, the identity (\ref{eq:rotWhit-id}) follows immediately from (\ref{eq:whit1-id}).
For (\ref{eq:rotWhit-x}), observe that setting $\omega:=1\in \mcH_0\Lambda^2(\R^2)$, we have that $\Phi(\kappa\omega)=\rot$.
Therefore, (\ref{eq:whit1-x}) implies that
\[\sum_{i<j}\left(\rot~\bv_i)\cdot \bv_j\right) \whit_{ij} = \rot~\bx.\]
Applying $\rot$ to both sides of the above equation completes the proof.
\end{proof}
\text{}\\

\begin{theorem}
\label{thm:pr-whit-ijk}
Let $\elt$ be a convex polyhedron in $\R^3$ with vertex set \textnormal{$\{\bv_i\}$} and an associated set of generalized barycentric coordinates $\{\lambda_i\}$.
Then
\begin{equation}
\label{eq:whit3-id}
\textnormal{$
\sum_{i<j<k}\whit_{ijk} \left((\bv_j-\bv_i)\times(\bv_k-\bv_i)\right)^T = \id,
$}
\end{equation}
and
\begin{equation}
\label{eq:whit3-x}
\textnormal{$
\sum_{i<j<k}(\bv_i\cdot(\bv_j\times\bv_k))\whit_{ijk}=\bx.
$}
\end{equation}
Thus, $\ds\textnormal{span}\left\{\whit_{ijk}\;:\;\textnormal{$\bv_i,\bv_j,\bv_k\in\elt$}\right\}\supseteq \cP_1^-\Lambda^2(\elt).$
\end{theorem}
\begin{proof}
We adopt the shorthand notations
\[\xi_{ijk} := \lambda_i\nabla\lambda_j\times\nabla\lambda_k,
  \quad \bz_{ijk}:= (\bv_j-\bv_i)\times(\bv_k-\bv_i), 
  \quad \bv_{ijk}:= \bv_i\cdot(\bv_j\times\bv_k).\]
For (\ref{eq:whit3-id}), we re-write (\ref{eq:lnlxnl-id}) as
\[
\sum_{i,j,k}\xi_{ijk}{\bz_{ijk}}^T= 2\id.
\]
Observe that $\xi_{ijk}{\bz_{ijk}}^T=(-\xi_{ikj})(-\bz_{ikj})^T=\xi_{ikj}{\bz_{ikj}}^T$ and $\bz_{ijk}=0$ if $i$, $j$, $k$ are not distinct.
Thus,
\begin{align*}
2\id
  & = \sum_{\substack{i<j<k \\ k<i<j \\ j<k<i}}\xi_{ijk}{\bz_{ijk}}^T + \sum_{\substack{i<k<j \\ k<j<i \\ j<i<k}}\xi_{ikj}{\bz_{ikj}}^T.
\end{align*}  
The two summations have different labels for the indices but are otherwise identical.
Therefore,
\begin{align*}
\id 
  & = \sum_{i<j<k}\xi_{ijk}{\bz_{ijk}}^T + \sum_{k<i<j}\xi_{ijk}{\bz_{ijk}}^T + \sum_{j<k<i}\xi_{ijk}{\bz_{ijk}}^T \\
  & = \sum_{i<j<k}\xi_{ijk}{\bz_{ijk}}^T + \xi_{jki}{\bz_{jki}}^T + \xi_{kij}{\bz_{kij}}^T \\
  & = \sum_{i<j<k}(\xi_{ijk}+\xi_{jki}+\xi_{kij}){\bz_{ijk}}^T  \\
  & = \sum_{i<j<k}\whit_{ijk}\left((\bv_j-\bv_i)\times(\bv_k-\bv_i)\right)^T.
\end{align*}
For (\ref{eq:whit3-x}), we take $\A$ as the identity, and re-write (\ref{eq:lnlxnl-Ax}) as
\[
\sum_{i,j,k}\bv_{ijk}\xi_{ijk} =2\bx.
\]
Observe that $\bv_{ijk}\xi_{ijk}=(-\bv_{ikj})(-\xi_{ikj})=\bv_{ikj}\xi_{ikj}$ and $\bv_{ijk}=0$ if $i$, $j$, $k$ are not distinct.
Thus,
\begin{align*}
2\bx 
  & = \sum_{\substack{i<j<k \\ k<i<j \\ j<k<i}}\bv_{ijk}\xi_{ijk} + \sum_{\substack{i<k<j \\ k<j<i \\ j<i<k}}\bv_{ikj}\xi_{ikj} .
\end{align*}  
The rest of the argument follows similarly, yielding
\begin{align*}
\bx  
  & = \sum_{i<j<k}\bv_{ijk}\xi_{ijk} + \sum_{k<i<j}\bv_{ijk}\xi_{ijk} + \sum_{j<k<i}\bv_{ijk}\xi_{ijk} \\
  & = \sum_{i<j<k}(\bv_i\cdot(\bv_j\times\bv_k))\whit_{ijk}.
\end{align*}
Note that $\mcH_{0}\Lambda^{3}(\elt)$ is generated by the volume form $\eta=dxdydz$ and that $\kappa\eta$ has vector proxy $\bx$.
Thus, by (\ref{eq:prminus-decomp}), we have shown that the span of the $\whit_{ijk}$ functions contains the vector proxy of any element of $\cP_1^-\Lambda^2(\elt)$.
\end{proof}
\text{}\\
\begin{remark}
\label{rmk:topdim}
\em
There are some additional constructions in this same vein that could be considered.
On a polygon in $\R^2$, we can define $\whit_{ijk}$ in the same way as (\ref{eq:whit-tri-def}), interpreting $\times$ as the two dimensional cross product.
Likewise, on a polyhedron in $\R^3$, we can define $\whit_{ijk\ell}$ according to formula (\ref{eq:whit-def}), yielding functions that are summations of terms like $\lambda_i(\nabla\lambda_j\cdot(\nlxnl k\ell)$.
These constructions will yield the expected polynomial reproduction results, yet they are not of practical interest in finite element contexts, as we will see in the next section.
\end{remark}

\section{Polygonal and Polyhedral Finite Element Families}
\label{sec:polyg-fams}

Let $\mesh$ be a mesh of convex $n$-dimensional polytopes $\{\elt\}$ in $\R^n$ with $n=2$ or $3$ and assign some ordering $\textnormal{$\bv_1,\ldots,\bv_p$}$ to all the vertices in the mesh.
Fix an associated set of generalized barycentric hat functions $\textnormal{$\hat\lambda_1,\ldots,\hat\lambda_p$}$ as in Section~\ref{sec:global-cnty}.
In Table~\ref{tab:comp-bases}, we list all the types of scalar-valued and vector-valued functions that we have defined this setting.
When used over all elements in a mesh of polygons or polyhedra, these functions have global continuity and polynomial reproduction properties as indicated in Table~\ref{tab:results-summary}.

These two properties -- global continuity and polynomial reproduction -- are essential and intertwined necessities in the construction of $H\Lambda^k$-conforming finite element methods on \textit{any} type of domain mesh.
Global continuity of type $H\Lambda^k$ ensures that the piecewise-defined approximate solution is an element of the function space $H\Lambda^k$ in which a solution is sought.
Polynomial reproduction of type $\cP_1\Lambda^k$ or $\cP_1^-\Lambda^k$ ensures that the error between the true solution and the approximate solution decays linearly with respect to the maximum diameter of a mesh element, as measured in $H\Lambda^k$ norm.
On meshes of simplicial elements, the basis functions listed in Table~\ref{tab:comp-bases} are known and often used as local bases for the corresponding classical finite element spaces listed in Table~\ref{tab:fe-spaces}, meaning our approach recapitulates known methods on simplicial meshes.

\begin{table}[t]
\centering
\begin{tabular}{c|c|c|c|c|c}
n & k & space & \# construction & \# boundary & \# polynomial \\ \hline 
 &&&&&\\[-2mm]
2 & 0 & $\cP_1 \Lambda^0(\elt)$   & $v$             & $v$    & $3$\\
  &  & $\cP_1^- \Lambda^0(\elt)$   & $v$             & $v$    & $3$\\
 &&&&&\\[-2mm]
  & 1 & $\cP_1 \Lambda^1(\elt)$   & $\displaystyle 2{v \choose 2}$        & $2e$   & $6$\\
   &&&&&\\[-2mm]
  &   & $\cP_1^- \Lambda^1(\elt)$ & $\displaystyle {v \choose 2}$ & $e$    & $3$\\
   &&&&&\\[-2mm]
  & 2 & $\cP_1 \Lambda^2(\elt)$   & $\displaystyle 3{v \choose 3}$ & $0$    & $3$\\
   &&&&&\\[-2mm]
  &   & $\cP_1^- \Lambda^2(\elt)$ & $\displaystyle {v \choose 3}$ & $0$    & $1$\\ [3mm]
  \hline
  \end{tabular}\\
\caption{Dimension counts relevant to serendipity-style reductions in basis size are shown. 
Here, $v$ and $e$ denote the number of vertices and edges in the polygonal element $\elt$.
The column `\# construction' gives the number of basis functions we define
(cf.\ Table~\ref{tab:comp-bases}), `\# boundary' gives the number of basis functions related to inter-element continuity, and `\# polynomial' gives the dimension of the contained space of polynomial differential forms.}
\label{tab:sizes-2d}
\end{table}

\paragraph*{\underline{Relation to quadrilateral and serendipity elements}}
Consider the scalar bi-quadratic element on \textit{rectangles}, which has nine degrees of freedom: one associated to each vertex, one to each edge midpoint, and one to the center of the square.
It has long been known that the `serendipity' element, which has only the eight degrees of freedom associated to the vertices and edge midpoints of the rectangle, is also an $H^1$-conforming, quadratic order method.  
In this case, polynomial reproduction requires the containment of $\cP_2\Lambda^0(\elt)$ in the span of the basis functions, meaning at least six functions are required per element $\elt\in\mesh$.
To ensure global continuity of $H^1$, however, the method must agree `up to quadratics' on each edge, which necessitates the eight degrees of freedom associated to the boundary.
Therefore, the serendipity space associated to the scalar bi-quadratic element on a rectangle has dimension eight.

In a previous paper~\cite{RGB2011a}, we generalized this `serendipity' reduction to $\cP_2\Lambda^0(\mesh)$ where $\mesh$ is a mesh of strictly convex polygons in $\R^2$.
For a simple polygon with $n$ vertices (and thus $n$ edges), polynomial reproduction still only requires $6$ basis functions, while global continuity of $H^1$ still requires reproduction of quadratics on edges, leading to a total of $2n$ basis functions required per element $\elt\in\mesh$.
Given a convex polygon, our approach takes the $n+{n \choose 2}$ pairwise products of all the $\lambda_i$ functions and forms explicit linear combinations to yield a set of $2n$ basis functions with the required global $H^1$ continuity and polynomial reproduction properties.

\paragraph*{\underline{Reduction of basis size}} A similar reduction procedure can be applied to the polygonal and polyhedral spaces described in Table~\ref{tab:comp-bases}.
A key observation is that the continuity results of Theorems \ref{thm:hcurl-conf} and \ref{thm:hdiv-conf} only rely on the agreement of basis functions whose indices are of vertices on a shared boundary edge (in 2D) or face (in 3D). 
For example, if vertices $\bv_i$ and $\bv_j$ form the edge of a polygon in a 2D mesh, $H(\Curl,\mesh)$ continuity across the edge comes from identical tangential contributions in the $\lnl ij$ and $\lnl ji$ functions from either element containing this edge and zero tangential contributions from all other basis functions.
Thus, basis functions whose indices do not belong to a single polygon edge (in 2D) or polyhedral face (in 3D) do not contribute to inter-element continuity, allowing the basis size to be reduced.

\begin{table}[ht]
\centering
\begin{tabular}{c|c|c|c|c|c}
n & k & space & \# construction & \# boundary & \# polynomial \\ \hline 
 &&&&&\\[-2mm]
3 & 0 & $\cP_1 \Lambda^0(\elt)$       & $v$             & $v$    & $4$\\
   &   & $\cP_1^- \Lambda^0(\elt)$   & $v$             & $v$    & $4$\\
   & 1 & $\cP_1 \Lambda^1(\elt)$   & $\displaystyle 2 {v \choose 2}$        & $\displaystyle \left(\sum_{a=1}^f v_a(v_a-1)\right) - 2e$   & $12$\\
   &   & $\cP_1^- \Lambda^1(\elt)$ & $\displaystyle {v \choose 2}$ & $\displaystyle \left(\sum_{a=1}^f {v_a\choose 2}\right) - e$    & $6$\\
   &&&&&\\[-3mm]
   & 2 & $\cP_1 \Lambda^2(\elt)$   & $\displaystyle 3{v \choose 3}$ & $\displaystyle \sum_{a=1}^f \frac{v_a(v_a-1)(v_a-2)}{2}$   & $12$\\
   &&&&&\\[-3mm]
   &   & $\cP_1^- \Lambda^2(\elt)$ & $\displaystyle {v \choose 3}$ & $\displaystyle \sum_{a=1}^f {v \choose 3}$    & $4$\\
   & 3 & $\cP_1 \Lambda^3(\elt)$   & $\displaystyle 4{v \choose 4}$ & $0$ & $4$\\
   &&&&&\\[-2mm]
   &   & $\cP_1^- \Lambda^3(\elt)$ & $\displaystyle {v \choose 4}$ & $0$    & $1$\\[3mm]
    \hline
\end{tabular}\\
\caption{The $n=3$ version of Table~\ref{tab:sizes-2d}.
Here, $f$ denotes the number of faces on a polyhedral element $\elt$ and $v_a$ denotes the number of vertices on a particular face $\face_a$.
The entries of the `\# boundary' column are determined by counting functions associated to each face of the polyhedron and, in the $k=1$ cases, accounting for double-counting by subtraction.
}
\label{tab:sizes-3d}
\end{table}

To quantify the extent to which the bases we have defined could be reduced without affecting the global continuity properties, we count the number of functions associated with codimension 1 faces for each space considered.
For a polygon in 2D, the results are summarized in Table~\ref{tab:sizes-2d}.
The $k=0$ case is optimal in the sense that every basis function $\lambda_i$ contributes to the $H^1$-continuity in some way, meaning no basis reduction is available.
In the $k=1$ cases, the number of basis functions we construct is quadratic in the number of vertices, $v$, of the polygon, but the number associated with the boundary is only linear in the number of edges, $e$.
Since $e=v$ for a simple polygon, this suggests a basis reduction procedure would be both relevant and useful; the description of such a reduction will be the focus of a future work.
In the $k=2$ cases, our procedure constructs $O(v^3)$ basis functions but no inter-element continuity is required; in these cases, a discontinuous Galerkin or other type of finite element method would be more practical.

For a polyhedron $\elt$ in 3D, the results are summarized in Table~\ref{tab:sizes-3d}.
As in 2D, the basis for the $k=0$ case cannot be reduced while the bases for the $k=n$ cases would not be practical for implementation since no inter-element continuity is required.
In the $k=1$ cases, the number of basis functions we construct is again quadratic in $v$, while the number of basis functions required for continuity can be reduced for non-simplicial polyhedra.
For instance, if $\elt$ is a hexahedron, our construction for $\cP_1\Lambda^1$  gives 56 functions but only 48 are relevant to continuity; in the $\cP_1^-\Lambda^1$ case, we construct 28 functions but only 20 are relevant to continuity.
In the $k=2$ cases, a similar reduction is possible for non-simplicial polyhedra.
Again in the case of a hexahedron, we construct 168 functions for $\cP_1\Lambda^1$ and 56 functions for $\cP_1^-\Lambda^1$, but the elements require only 72 and 24 functions, respectively, for inter-element continuity.

\paragraph*{\underline{Current and future directions}}
It remains to discover additional properties of Whitney-like basis functions built from generalized barycentric coordinates and their use in finite element methods.
In the time since this manuscript first appeared online, Chen and Wang~\cite{CW2015} have presented an approach for constructing `minimal dimension' local basis sets based on the results of this paper.
Their theoretical and numerical results indicate that minimal spaces can, indeed, be constructed using the methods presented here with expected rates of convergence on certain classes of polygons and polyhedra.
We expect that the ideas introduced here will continue to influence the rapidly expanding use of polytopal finite element methods in scientific and engineering applications.\\

\noindent
\textbf{Acknowledgements.}
The authors would like to thank the anonymous referees for their helpful suggestions to improve the paper.
AG was supported in part by NSF Award 1522289 and a J.\ Tinsley Oden Fellowship.
CB was supported by was supported in part by a grant from NIH (R01-GM117594) and contract (BD-4485) from Sandia National Labs.
Sean Stephens helped produce the figure.

\bibliographystyle{abbrv}
\bibliography{references}
\end{document}